%% file: impl_needed_art.tex
\documentclass[11pt]{article}
\usepackage[a4paper]{geometry}
\usepackage{amsfonts,amsmath,amssymb,bm,graphicx,hyperref,xcolor}

\title{Are implicit schemes needed in adaptive time stepping\\for PDEs?%
\thanks{To the memory of Howard Elman.}}
\author{Mikhail A. Botchev\thanks{Keldysh Institute of Applied Mathematics of Russian Academy of Sciences, Moscow, 125047 Russia, \url{botchev@kiam.ru}}}

\begin{document}
\maketitle

\begin{abstract}
\input{abs.txt}


\end{abstract}

\input{impl_needed_cont_rev}

\bibliography{my_bib,applic,matfun}
\bibliographystyle{ieeetr}

\end{document}

%% file: abs.txt
% Are implicit schemes needed in adaptive time accurate PDE integration?

Implicit time integration is the standard way to overcome the 
tight CFL time step restriction when solving nonstationary PDEs numerically.
However, within an adaptive time accurate integration, the time step size 
selection is governed not only by stability but also by accuracy requirements.
Therefore, an eligible question to ask is: is (unconditionally) stable 
implicit time stepping really needed if the time step is restricted by 
accuracy considerations anyway?
And how much do we then lose (in terms of the time step size) when switching 
from the expensive implicit time stepping to stabilized explicit schemes?
An answer to this question is also important because such a switch is
often done when solving large scale nonstationary PDEs on 
%\textcolor{blue}{
massively
%} 
parallel computers.

This paper attempts to answer these questions by numerical tests 
for two particular adaptive time accurate settings:
(i)~nonstationary convection-diffusion equations 
and 
(ii)~nonlinear heat equations.

%% file: impl_needed_cont_rev.tex
\maketitle

\newcommand{\Cc}{\mathbb{C}}
\newcommand{\eps}{\epsilon}
\newcommand{\epst}{\tilde{\varepsilon}}
\newcommand{\et}{\tilde{e}}
\newcommand{\geqs}{\geqslant}
\newcommand{\lambmax}{\lambda_{\max}}
\newcommand{\leqs}{\leqslant}
\newcommand{\Rr}{\mathbb{R}}
\newcommand{\vareps}{\varepsilon}
\newcommand{\yh}{\hat{y}}
\newcommand{\yt}{\tilde{y}}

\newtheorem{proposition}{Proposition}
\newtheorem{remark}{Remark}

\section{Introduction}
Implicit time integration schemes have been widely employed in 
computational fluid dynamics (CFD) simulations~\cite{bookWesseling}.  
They proved to 
%\textcolor{blue}{
be
%} 
a useful means for solving stiff problems, 
for which 
explicit schemes would have been inefficient due to the prohibitively
tight stability restrictions on the time step size~\cite{HundsdorferVerwer:book}.  
However, within 
an adaptive implicit time integration setting, the time step size is 
restricted not only by stability but also by the accuracy requirements.  
Therefore, an eligible answer to ask is whether the time step size used 
by an implicit scheme is within the reach of stabilized explicit schemes.
A possible switch to cheap (per time step) explicit schemes will then
%\textcolor{blue}{
not only
%} 
save computational work 
%\textcolor{blue}{
but also
%} 
make an implementation on a parallel or hybrid computer much easier.    

%% %Correction to be hopefully made in the proofs: !!!
In general implicit time integration schemes 
are not trivial to implement on massively parallel platforms,
as large nonlinear or linear systems have to be solved in parallel.
However,   
last decades significant advances have been made in developing parallel 
linear solvers, 
% end of corrections for the proofs
in particular, in 
domain decomposition methods~\cite{Tang_ea2020_DDM}, 
parallel 
preconditioners~\cite{Kaporin2012,AnztHuckle_ea2018,PearsonPestana2020} and
parallel sparse direct 
methods~\cite{BollhoeferSchenkJanalikHammGullapalli2020,LiLiu2026}.
Additional parallelism can also, to some extend, be gained by
employing parallel-across-time 
schemes~\cite{Burrage1997,vdHouwen_deSwart1997,Gander2015,Kooij_ea2017,Kooij_ea2018}.
Nevertheless, 
we note that in 
%\textcolor{blue}{
massively
%} 
parallel implementations 
of large scale PDE simulations,
explicit stabilized schemes are often the schemes of 
choice~\cite{ChetverushkinOlkhovskaya2020,Balashov_ea2020,RichterSchopsClemens2017gpu,Zhukov_ea2020}.  
Hence, answering the question how large the gap is between
implicit and stabilized explicit time stepping, allows to evaluate
how much we actually lose, in terms of the time step sizes, 
with respect to implicit time stepping. 
 
This note attempts to give an answer to these questions by a numerical
study in two particular adaptive time integration settings.  
The first one is the unsteady convection--diffusion problem discretized
in space by finite elements, provided by the IFISS
(Incompressible Flow \&{} Iterative Solver Software) software 
library~\cite{ElmanSilvesterWathen:book2,IFISS_sirev}.
The second setting is a nonlinear heat equation, discretized in space
by standard second order finite differences, as considered
in~\cite[App.~1, Ch.~2.11]{TikhonovSamarskii_UMF} and~\cite{BotchevZhukov2024}.
We choose these test problems because they seem to present two major 
important classes of parabolic equations often arising in CFD simulations.
Indeed, though being relatively simple, convection-diffusion problems,
with space-discretized operators being nonsymmetric matrices, 
typically form a useful prototype for the Navier-Stokes equation and
``can exhibit multiple time scales''~\cite{GreshoGriffithsSilvester2008}.
On the other hand, the nonlinear heat equations,
where space-discretized operators are symmetric matrices with usually
a stiff spectrum, are often encountered in various applications such as, 
e.g., 
three-temperature gas dynamics simulations~\cite{ZabrodinProkopov1998}
numerical models of non-steady gas dynamics processes~\cite{Feodoritova_ea2021}, and high-temperature superconductivity physics~\cite{Feodoritova_ea2022}.  

This paper is structured as follows.
The next section discusses the experimental setting and assumptions
made in this study.
In Section~\ref{s:cd}, adaptive time accurate integration
of non-steady convection-diffusion problems is considered.
Section~\ref{s:heat} is devoted to adaptive time accurate integration
of nonlinear heat equation.
Conclusions are formulated in Section~\ref{s:concl}.

\section{Experiments setting and assumptions}
\label{s:req}
As outlined in the introduction, one of the aims of this paper 
is to evaluate performance of stabilized explicit schemes,
which are attractive in 
%\textcolor{blue}{
massively
%} 
parallel computing environments.
Therefore, in the tests presented here we restrict our attention
to explicit schemes with ``light-weighted'' iterative solvers only.
More specifically, this means that we do not consider applying possible
preconditioning methods, except that, when appropriate, the simplest
diagonal (Jacobi) preconditioner is used.
Modern preconditioning techniques can be very efficient and lead 
to a fast convergence of iterative 
methods~\cite{vdVBook,SaadBook2003,Benzi02prec,Wathen2015prec}.  
However, in terms of
parallel implementation, a powerful preconditioner can be just
as difficult to implement as a (sparse) direct solver. 
Indeed, even simple 
triangular preconditioners such as the standard incomplete LU (ILU) 
factorizations or SSOR preconditioners, do not allow a straightforward
parallelization, see, e.g.,~\cite{Henk:HPP,NLA-SIAM,vdVBook,SaadBook2003}.
Furthermore, we restrict ourselves to short-term recurrence methods
(in particular, in Section~\ref{s:cd} we use 
the BiCGstab method~\cite{BiCGStab,vdVBook})
and Krylov subspace methods for evaluating matrix exponential actions
with limited memory requirements (Section~\ref{s:heat}).
The restricted memory usage in the latter methods is reached by applying 
an efficient restarting~\cite{BKT21}. 
If the coefficient matrix of the ODE system is symmetric (as in the
setting of Section~\ref{s:heat}) then Chebyshev polynomials
can be used for an efficient explicit time stepping.  
An example of such schemes is the LIM time integration 
method~\cite{LokLokDAN,Zhukov2011}
which is included in numerical comparisons of Section~\ref{s:heat}.

It should be emphasized that such ``explicit'' iterative techniques
do have intrinsic limitations.  For instance, 
if a one-space-dimensional parabolic problem is solved and the
space discretization is built on a standard three-point finite difference stencil,
then, to propagate information through the whole space domain,
it is clear that approximately $n$  matrix-vector multiplications
are needed, with $n$ being the number of space grid points.
This fundamental limitation was observed already in 1952 by
Gel'fand and Lokutsievskii, who, in an appendix~\cite{GelfandLokutsievskii1962} 
to the pioneering book~\cite{GodunRyab1964}, proved that
approximately $\sqrt{2\tau/h^2}$ matrix-vector multiplications
are required to have stability for the time step size~$\tau$ and
space grid size~$h=1/n$.

All the numerical tests described in this paper are carried out
in the GNU Octave programming language, see~\url{https://octave.org/}.

\begin{figure}
\includegraphics[width=0.51\linewidth]{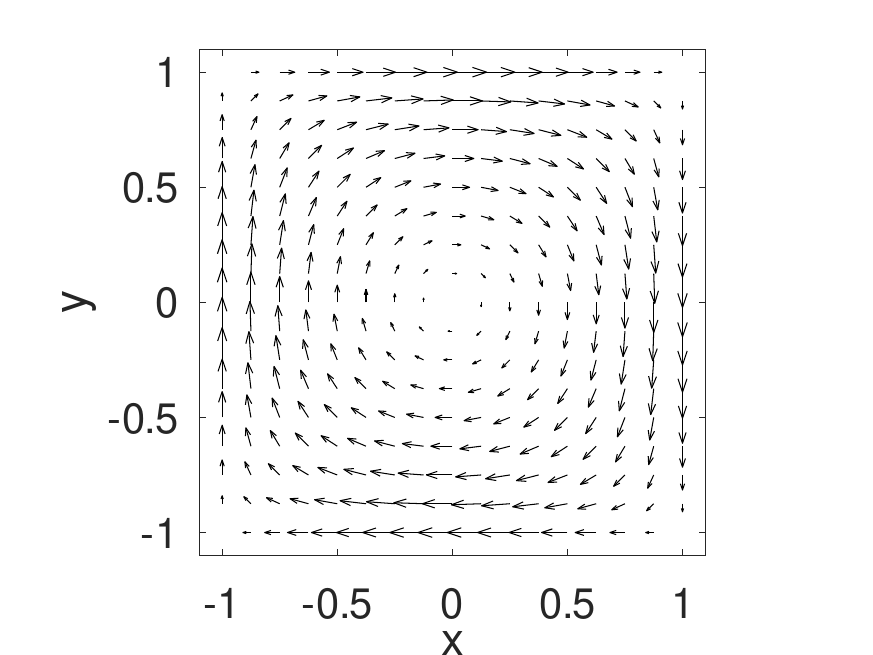}%
\includegraphics[width=0.48\linewidth]{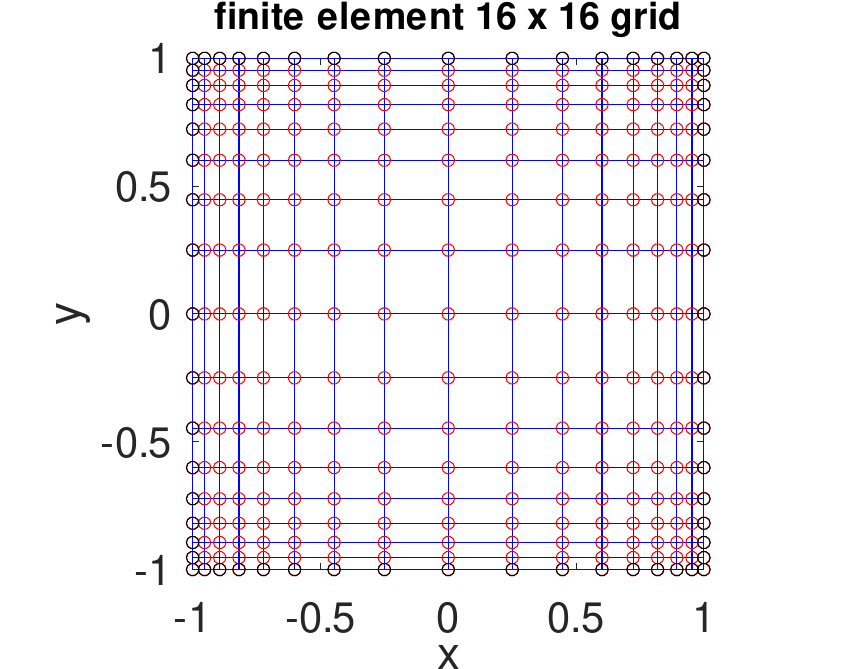}
\caption{Left: a ``recirculating wind'' velocity vector field~\eqref{wind}
from the IFISS test problem.  
Right: the stretched $16\times 16$ IFISS finite element grid.}  
% adaptively chosen time step size values $\tau$ for the $256\times 256$ grid
% and viscosity $\nu=1/200$.}
\label{f:vel_grid}
\end{figure}

\section{Nonsteady convection--diffusion problem}
\label{s:cd}
\begin{figure}
\includegraphics[width=\linewidth]{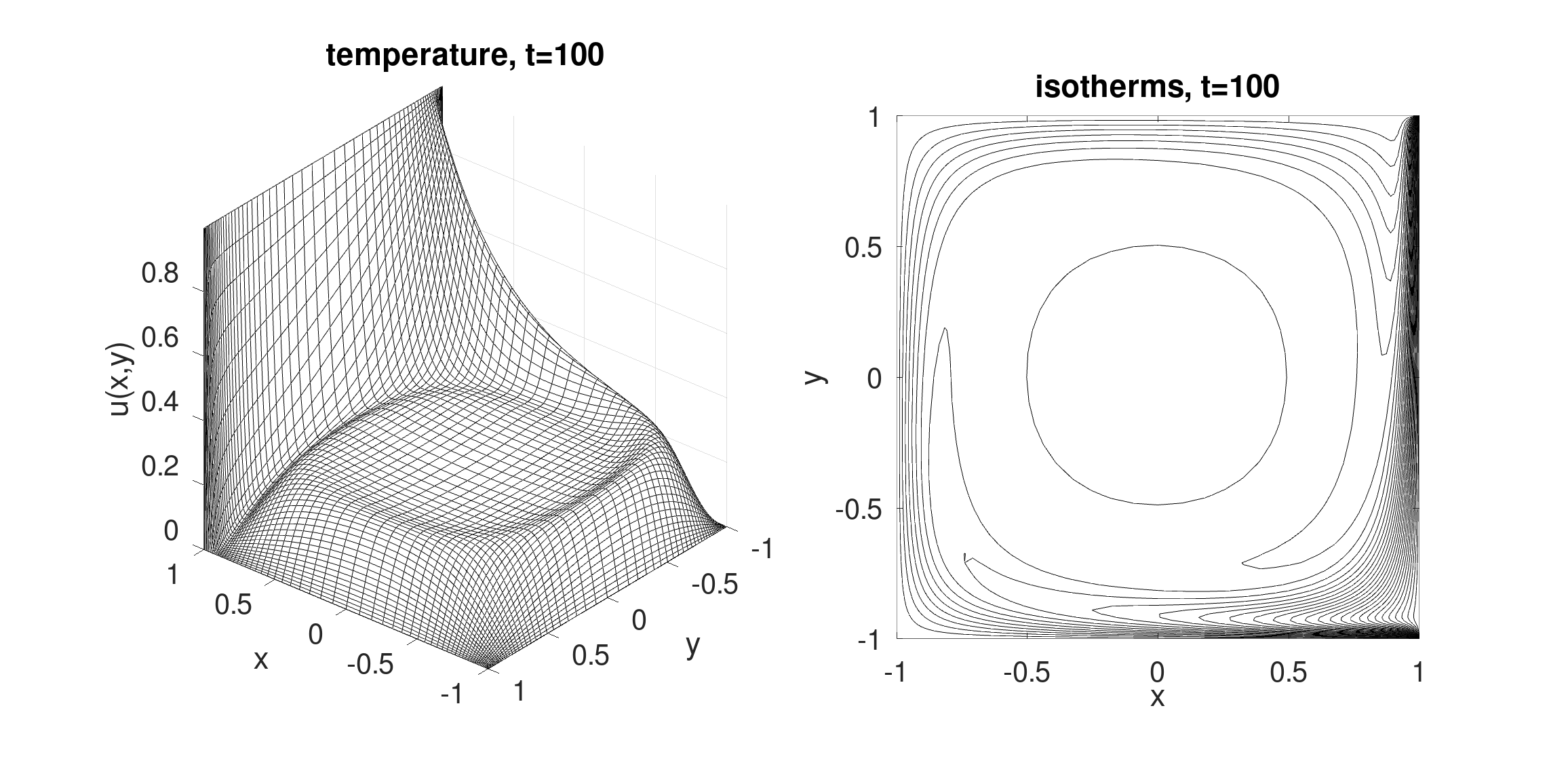}
\caption{Solution of the IFISS convection-diffusion test problem on the stretched $128\times 128$ mesh for viscosity $\nu=1/200$ at time $t=100$}
\label{f:cd_sol}
\end{figure}

This test is carried out with the help of
the IFISS (Incompressible Flow \&{} Iterative Solver Software) software 
library~\cite{ElmanSilvesterWathen:book2,IFISS_sirev},
freely available for download at \url{http://www.manchester.ac.uk/ifiss}.
Initial-boundary value problem being solved,
for unknown temperature distribution $u=u(x,y,t)$, reads
\begin{equation}
\label{cd}
\begin{gathered}
\frac{\partial u}{\partial t} - \nu\Delta u + \bm{w}\cdot\nabla u = 0,
\quad (x,y)\in\Omega=[-1,1]\times[-1,1],
\\
u(x,y,0)=0, \quad u(-1,y,t)=u(x,\pm y,t)=0, \quad u(1,y,t)= 1,
\end{gathered}
\end{equation}
where $\nu$ is viscosity and $\bm{w}=\bm{w}(x,y)$ is 
velocity field vector function
\begin{equation}
\label{wind}
\bm{w}(x,y) = 
\begin{pmatrix}
2y(1 - x^2 ), & -2x(1 - y^2 )
\end{pmatrix},
\end{equation}
whose vector field plot is given in Figure~\ref{f:vel_grid}.
Thus, this problem models temperature distribution evolution,
for zero initial value and homogeneous Dirichlet boundary conditions,
with only the wall $x=1$ being kept hot.
%% %Correction to be hopefully made in the proofs: !!!
A typical solution at time~$t=100$ is shown % at the left plot 
in Figure~\ref{f:cd_sol}.

Taken into account that the large velocity values 
(and, as a consequence, significant values of the solution gradients)
near the domain boundary, it is sensible to use a nonuniform grid.
An example of the stretched grid used in the tests is shown at the 
right plot in Figure~\ref{f:vel_grid}.  
For the grids used in the tests, their stretch ratios are reported 
in Table~\ref{t:grid_nu_Pe}.   
The IFISS package discretizes problem~\eqref{cd} in space by the
bilinear quadrilateral $Q_1$ finite element approximation,
stabilized by the streamline diffusion method.  This means that some 
artificial diffusion is added in the streamline direction.  The amount
of the added diffusion depends on the element Peclet number, defined
for a finite element~$k$ as  
$\textit{Pe}=|\bm{w_k}|h_k/(2\nu)$, where $|\bm{w_k}|$ is the wind length
at the element centroid and $h_k$ is an estimate of the element size
in the wind direction, see~\cite{IFISSacm2007} for more details.
The maximal element Peclet numbers, for the grids and viscosity values
used in the tests, are given in Table~\ref{t:grid_nu_Pe}. 
                
\begin{table}
\caption{Grid sizes, stretch ratios, viscosity values used in the tests and 
corresponding maximal 
%% %Correction to be hopefully made in the proofs: !!!
element
Peclet numbers, nonsymmetry estimates
for the matrix $R=M+\frac\tau2A$, $\tau=\max_n\tau_n$,
and the spectral radii of the ODE system matrix $M^{-1}A$}
\label{t:grid_nu_Pe}
\begin{center}
\addtolength{\tabcolsep}{2ex}
\centerline{\begin{tabular}{cccccc}
\hline\hline
grid  & grid stretch & viscosity & max. element  & 
$\frac{\|R-R^T\|_1}{\|R+R^T\|_1}$ & 
$\rho(M^{-1}A)$
\\[-2   ex]
size  &  ratio       &  $\nu$    & Peclet number &                          
\\\hline
$256\times 256$ & 1.032 & 1/200 & 6.247 & {\tt2.182e-02} & {\tt2.493e+05}
\\  
$256\times 256$ & 1.032 & 1/50  & 1.562 & {\tt5.454e-03} & {\tt9.971e+05}
\\  
$512\times 512$ & 1.018 & 1/200 & 3.515 & {\tt7.104e-03} & {\tt2.414e+06}
\\
$512\times 512$ & 1.018 & 1/50  & 0.879 & {\tt1.776e-03} & {\tt9.654e+06}
\\\hline
\end{tabular}}
\end{center}
\end{table}

\begin{figure}
\includegraphics[width=0.49\linewidth]{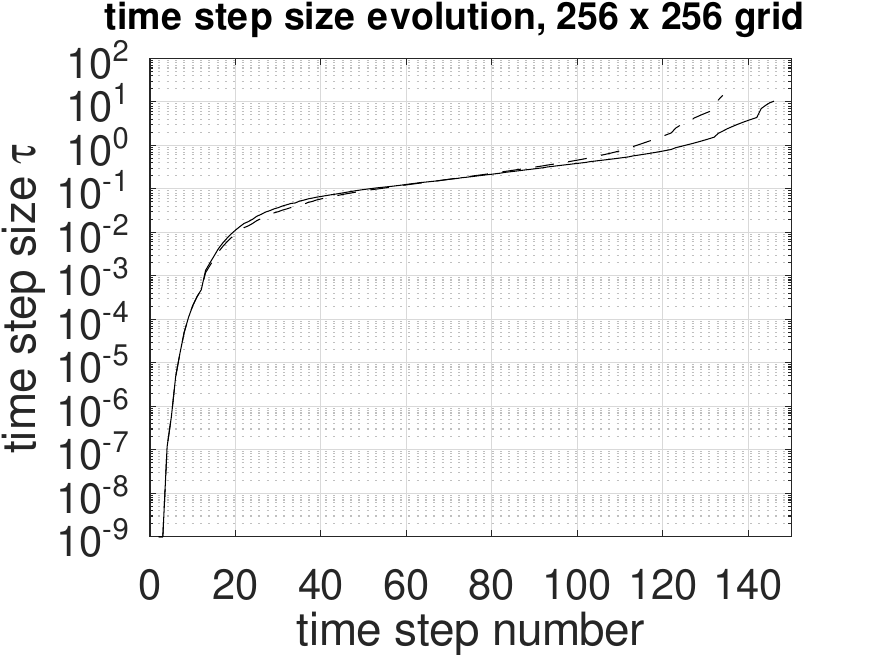}%
\includegraphics[width=0.49\linewidth]{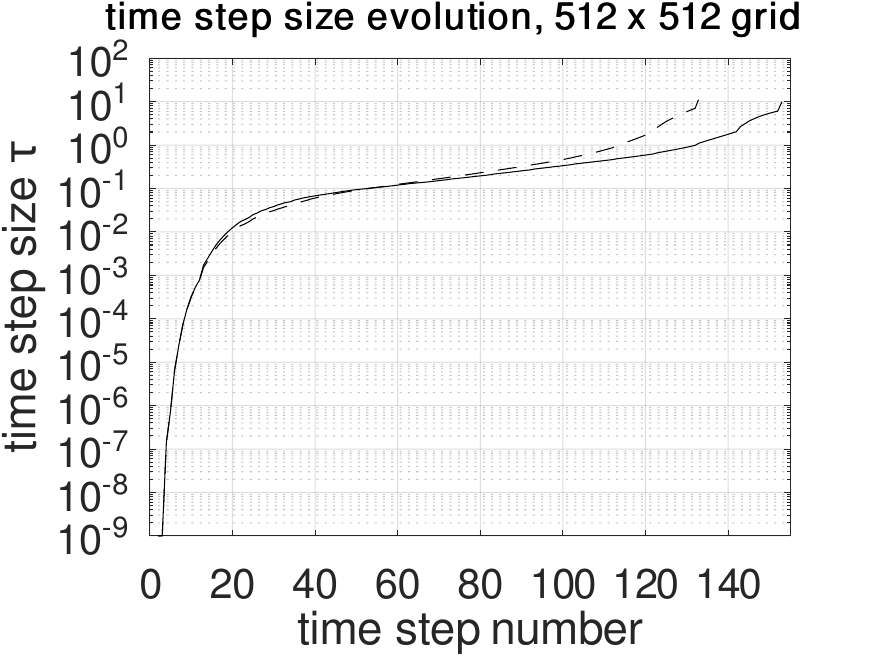}
\caption{Time step size evolution in the adaptive
implicit trapezoidal rule for the
$256\times 256$ (left plot) and $512\times 512$ (right plot) 
finite element grids, with viscosity values 
$\nu=1/200$ (solid line) and $\nu=1/50$ (dashed line).
The accuracy tolerance in the adaptive time stepping is
in all cases $10^{-4}$.}
\label{f:dt}
\end{figure}

The $Q_1$ stabilized finite element space approximation reduces the initial-value
boundary problem~\eqref{cd} to an initial value problem (IVP) 
\begin{equation}
\label{cd_ivp}
M y'(t) = -Ay(t) + g, \quad y(0)=0,
\end{equation}
where $y(t)$ is the unknown vector function containing the finite element
degrees of freedom, 
$M\in\Rr^{N\times N}$ is a symmetric positive definite finite element
mass matrix, and $A\ne A^T\in\Rr^{N\times N}$ is the stiffness matrix approximating 
the convection-diffusion operator, 
i.e., $Ay\approx - \nu\Delta u + \bm{w}\cdot\nabla u$.   
The IVP~\eqref{cd_ivp} is then integrated in time by an adaptive
implicit trapezoidal rule scheme 
\begin{equation}
\label{itr}
M\frac{y^{n+1}-y^n}{\tau_n} = -\frac12 ( Ay^n + Ay^{n+1}) + g,
\quad n=0,1,\dots,
\end{equation}
where $n$ is the time step number and $\tau_n>0$ is an adaptively chosen
time step size.  This adaptive time integration scheme,
where the time step size selection is based on a combination
with the second order Adams-Bashforth method, is discussed 
in detail in~\cite{GreshoGriffithsSilvester2008}.
At every time step~$n$ in~\eqref{itr} a linear system in~$y^{n+1}$
with the matrix
\begin{equation}
\label{R}
R = M + \frac{\tau_n}2 A
\end{equation}
has to be solved.
In Table~\ref{t:grid_nu_Pe}, we give the ratios 
$\|R-R^T\|_1/\|R+R^T\|_1$ showing the nonsymmetry degree of the matrices~$R$.
At the last column of the same table, the estimates for spectral
radius $\rho(M^{-1}A)$ are presented as a measure of the ODE stiffness
in~\eqref{cd_ivp}.
As we see in the table, the matrices $R$ can be considered as weakly nonsymmetric.
It should be noted that the nonsymmetry is small not because the mass matrix
$M$ is significant in norm relatively to $A$.  
Indeed, as can be seen in the last column of the same table, this is not the case.
Hence, the small nonsymmetry can be attributed to the streamline
finite element stabilization.
Examples of the time step size evolution in the course
of the adaptive time stepping are presented in Figure~\ref{f:dt}.
As we see, even with the accuracy tolerance $10^{-4}$ (the default
accuracy in IFISS for this test problem) the time step sizes reach values 
%% %Correction to be hopefully made in the proofs: !!!
$\tau_n\approx 10$, %10, 
for which we have $\tau_n\rho(M^{-1}A)\approx 10^7$ on the 
$512\times 512$ grid.  Therefore, this problem can be seen as highly
stiff.

There are two modifications we have made in the IFISS library
to carry out the numerical tests presented here.  
In IFISS, to impose the boundary conditions in the course of time stepping,
the linear system matrix $R$ (cf.~\eqref{R}) and the right hand side vector 
are modified in such a way that the next time step solution~$y^{n+1}$
gets the proper boundary values.  To achieve this, before the system
with $R$ is solved, the rows and columns 
in $R$ corresponding to boundary degrees of freedom are zeroed out,
while the diagonal entries are set to one.  The corresponding entries 
in the right hand side vector are set to the boundary values.
This procedure can significantly increase the norm of the matrix~$R$
which is harmful for convergence of iterative methods we apply.
To avoid this effect, we set the diagonal entries in~$R$ to
$\|R\|_1/10$ rather than to~one and adjust the entries in the
right hand side accordingly.

The second modification we have made in the IFISS library is that
the backslash operator, used to solve the system in~$y^{n+1}$
(i.e., the UMFPACK sparse direct solver~\cite{UMFPACK1,UMFPACK2}
in the backslash disguise), 
is replaced by an iterative solver.
Keeping in mind the parallelization prerequisites 
for the preconditioned iterative solvers discussed in Section~\ref{s:req},
we have compared several simple short-term recurrence iterative
methods, provided either by Octave or IFISS.
In particular, we have tested GMRES(10) (though short restarts
in GMRES can only loosely be considered as a ``short-term recurrence''),
BiCG, (TF)QMR, BiCGstab and BiCGstab(2) iterative solvers, 
all preconditioned from the right by the diagonal of $R$. 
These experiments revealed the diagonally preconditioned
BiCGstab as the most efficient, in terms
of required matrix-vector multiplications, solver.
Therefore, here we report the number matrix-vector multiplications
for this solver only.

\begin{figure}
\includegraphics[width=0.49\linewidth]{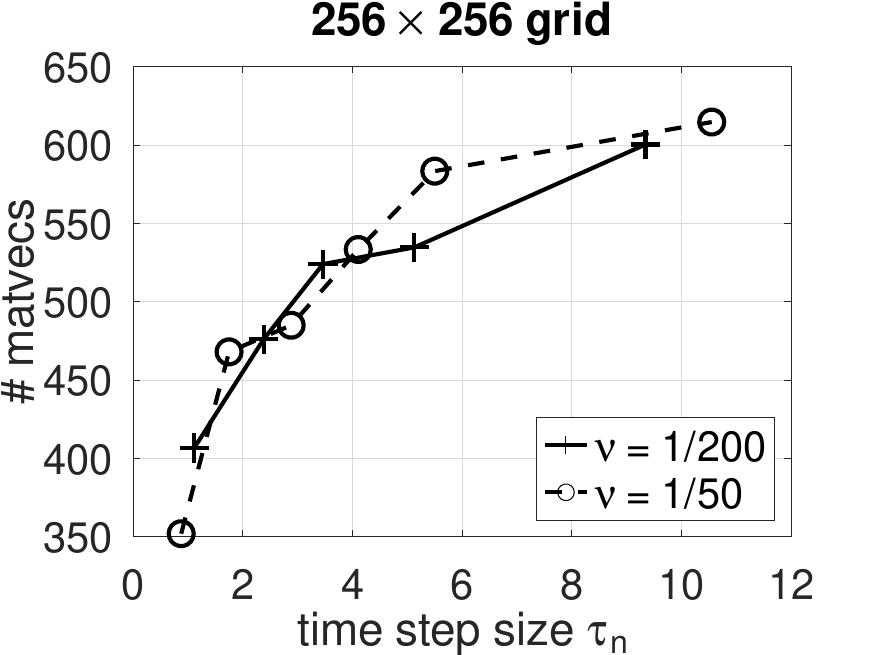}%
\includegraphics[width=0.49\linewidth]{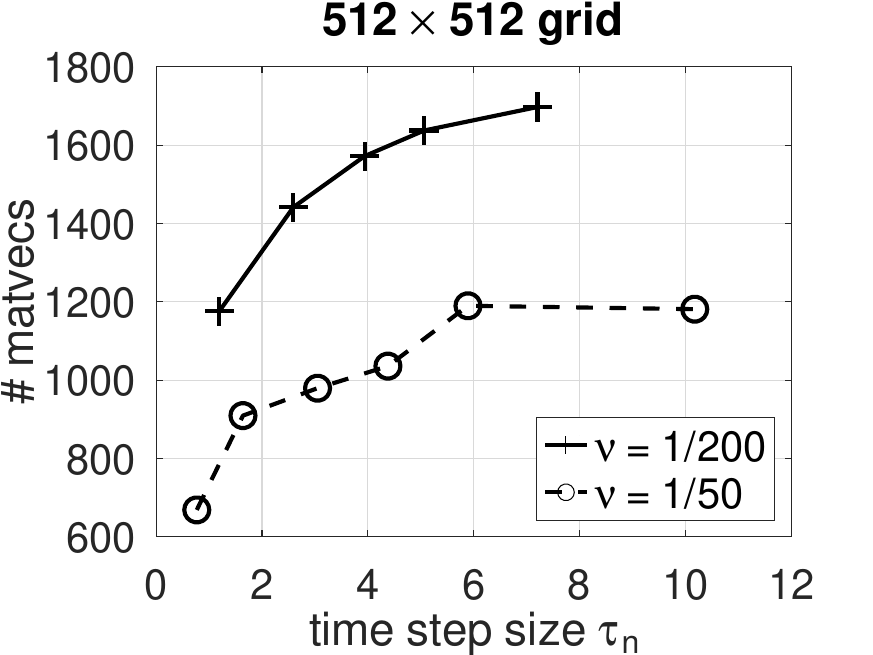}
\caption{Averaged numbers of the matrix-vector multiplications required 
by diagonally right-preconditioned BiCGstab for the
$256\times 256$ (left plot) and $512\times 512$ (right plot) 
finite element grids, with viscosity values 
$\nu=1/200$ (solid line) and $\nu=1/50$ (dashed line).
The stopping criterion tolerance in BiCGstab is set
to~$10^{-4}$, the accuracy tolerance in the
adaptive implicit time integrator.}
\label{f:cd}
\end{figure}

In the tests, the stopping criterion tolerance in the preconditioned BiCGstab 
is set to the default accuracy tolerance of the adaptive implicit trapezoidal rule 
scheme in IFISS, i.e., to $10^{-4}$. 
%% %Correction to be hopefully made in the proofs: !!!
This means that the iterations stop as soon as the residual is~$10^4$ times
smaller in norm than the linear system right hand side.
With this stopping criterion,
the relative difference in norm between the iterative and direct
solution vectors never exceeds $10^{-3}$.
The results of the tests are presented in Figure~\ref{f:cd}.
Since the BiCGstab typically yields an irregular convergence
behavior, the data in the figure are averaged, i.e., we plot
the mean value of three time step size values on the abscissa
axis and the mean of the corresponding matrix-vector multiplications
values on the ordinate axis. 

Several conclusions can be made from the data presented in Figure~\ref{f:cd}.
First, as we see, the convergence in the iterative solver depends on
the time step size rather weakly.
Second, we see that the solver performance heavily depends on the
grid size and, hence, the increased stiffness (see the last column
in Table~\ref{t:grid_nu_Pe}) plays the key role here. 
The role of the viscosity is much less evident.  On one hand, a larger
viscosity makes the system matrix~$R$ in~\eqref{R} more symmetric, which should
be favorable for the iterative solvers.  On the other hand, 
a larger viscosity value increases the diffusion with its
$h^{-2}$ dependence on the grid size and, hence, 
leads to a stiffer matrix (again, see the last column
in Table~\ref{t:grid_nu_Pe}).
This is apparently the reason why the viscosity values affect
the solver performance on the finer $512\times 512$ grid to a larger
extend than on the coarser grid.  
In addition to this subtle balance between the symmetry degree and the stiffness,
another aspect to consider is the streamline stabilization:
a less symmetric matrix may have obtained more stabilization.  

Since the presented here tests are carried out in Octave,
the CPU timings of the original adaptive implicit trapezoidal rule
and the scheme with the BiCGstab solver can not be fairly compared.
Indeed, the former uses the built-in precompiled direct solver which is quite
efficient for this 2D problem, whereas the iterative solvers in Octave 
are not precompiled.  It is clear that there is a gap in performance
between the implicit solver and the stabilized explicit solver with
BiCGstab iterations.  However, this gap turns out to be much smaller
in practice than could be expected: while the estimates of $\tau_n\rho(M^{-1}A)$,
showing the stiffness degree, reach approximately $10^8$ in this test,
the numbers of explicit iterations are of order $10^3$.

%\textcolor{blue}{
Although it is in general not easy to predict
how the implicit solver and the stabilized explicit solver with BiCGstab
iterations would perform in a parallel environment, 
several remarks can be made.
First, we note that nowadays there are modern sparse direct solvers
with ready-to-use distributed memory implementations available,
such as PARDISO (\url{https://www.pardiso-project.org/}) 
and MUMPS (\url{https://mumps-solver.org/}).
However, the parallelism in sparse direct solvers is often restricted,
especially for 3D problems~\cite{Kaporin2010},
and a reachable parallel efficiency is very problem dependent.  
For instance, as noted in~\cite{BollhoeferSchenkJanalikHammGullapalli2020},
``for a majority of simulations, scaling beyond eight cores is difficult''
in integrated circuit 
%% %Correction to be hopefully made in the proofs: !!!
modeling. % simulations.
%}

%\textcolor{blue}{
Furthermore, as already discussed above, ILU preconditioners
prevent parallelization in Krylov subspace iterative solvers.
Preconditioners allowing a better parallelism exist, including
parallel explicit and polynomial versions of ILU 
preconditioners~\cite{Henk:HPP,NLA-SIAM}, 
sparse approximate inverse preconditioners%
~\cite{Benzi02prec,Wathen2015prec,AnztHuckle_ea2018},
and domain decomposition based 
preconditioners~\cite{SaadBook2003,OlshTEE_book,Tang_ea2020_DDM}.
Sparse approximate inverse preconditioners 
``often fail to provide substantial convergence improvement''~\cite{AnztHuckle_ea2018} and the same can be stated about parallel explicit ILU and polynomial 
preconditioning~\cite{NLA-SIAM}.
Domain decomposition methods can certainly provide a good parallelism level 
but their implementation often requires 
significant efforts and tuning, for instance, in treating the arising
Schur complement systems properly~\cite{SaadBook2003,Tang_ea2020_DDM}.
Note that more advanced parallel versions of ILU can be seen as
generalized domain decomposition methods, see, e.g.,~\cite{Magolu_vdVorst2002}.
Finally, we mention that parallel linear solvers can be built based on the
multigrid approach, see, 
e.g.,~\cite{OlshTEE_book,FeodoritovaNovikovaZhukov2015}.
Summarizing, one can state that reliability of modern preconditioned
iterative solvers ``is still worse than that of direct solvers''~\cite{Kaporin2010},
whereas parallelism requirements in preconditioners
``imposes fairly heavy demands on their structure''~\cite{Kaporin2012}
and, hence, often lead to a performance deterioration.
%} 
      
%\textcolor{blue}{
We note that for 3D test problems the performance gap 
between the implicit and the stabilized explicit solvers
will likely be smaller for two reasons.
%}  
First, the direct sparse solvers
will be much less efficient than for 2D problems, see, e.g.,~\cite{Kaporin2010}.  
Second, the space
grid size $h$ in 3D problems 
%\textcolor{blue}{
can often be expected not as small as in 2D
problems, leading to better conditioned matrices.
Otherwise, for 3D problems to be treated on a fine space grid,
massively parallel computing platforms will likely have to be used,
where stabilized explicit iterative schemes may be more attractive than fully implicit schemes.
%}  

Taking into account the small degree of nonsymmetry in the matrix~$R$
in~\eqref{R}, 
we have also attempted to replace the BiCGstab iterations 
by a stable operator based on Chebyshev polynomials.
For instance, the LIM time step integration scheme, 
discussed in detail in the following section, provides such an operator.
In fact, the LIM scheme builds a polynomial in $\tau_nM^{-1}A$ which approximates 
the exact time propagator operator $\exp(-\tau_nM^{-1}A)$ for the modes
corresponding to the smallest eigenvalues of $M^{-1}A$, whereas
the other modes are just damped.
The number of the matrix-vector products 
per time step %% %Correction to be hopefully made in the proofs: !!!
in the LIM scheme
is roughly proportional to $\sqrt{\tau_n\rho(M^{-1}A)}$ 
and we actually see  that the same holds in this test for
BiCGstab: the switch to a finer $512\times 512$ grid means
an increase by a factor $\approx 10$ in $\rho(M^{-1}A)$ and
by a factor of $\approx 3$ in the matrix-vector products.
However, it turns out that the LIM scheme would need approximately $11\,000$ 
matrix-vector products for $\tau_n\approx 10$ and viscosity
$\nu=1/50$ on the $512\times 512$ grid.  This is a way more
than $\approx 1700$ matrix-vector products in BiCGstab.
This significant difference in the number of required matrix-vector
products can be attributed to the ability of the Krylov subspace
solvers to adapt to the discrete structure of the matrix spectrum
and to the initial (right-hand side) vector.
 
%% Another observation is that Krylov subspace iteration methods perform
%% here a way better than a stabilized explicit scheme based on Chebyshev
%% polynomials would perfrom (assuming that the matrix would be symmetric or
%% stabilized explicit scheme would be applied to the symmetric part of the 
%% matrix).  Indeed, ... give estimates of p(LIM)...

\section{Nonlinear heat equation}
\label{s:heat}
\begin{figure}
\includegraphics[width=0.49\linewidth]{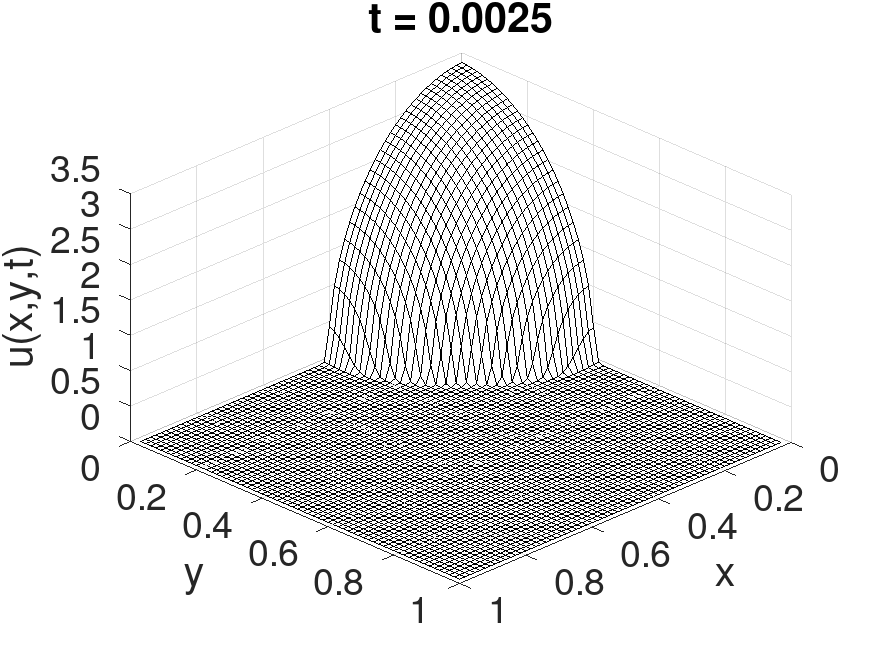}%
\includegraphics[width=0.49\linewidth]{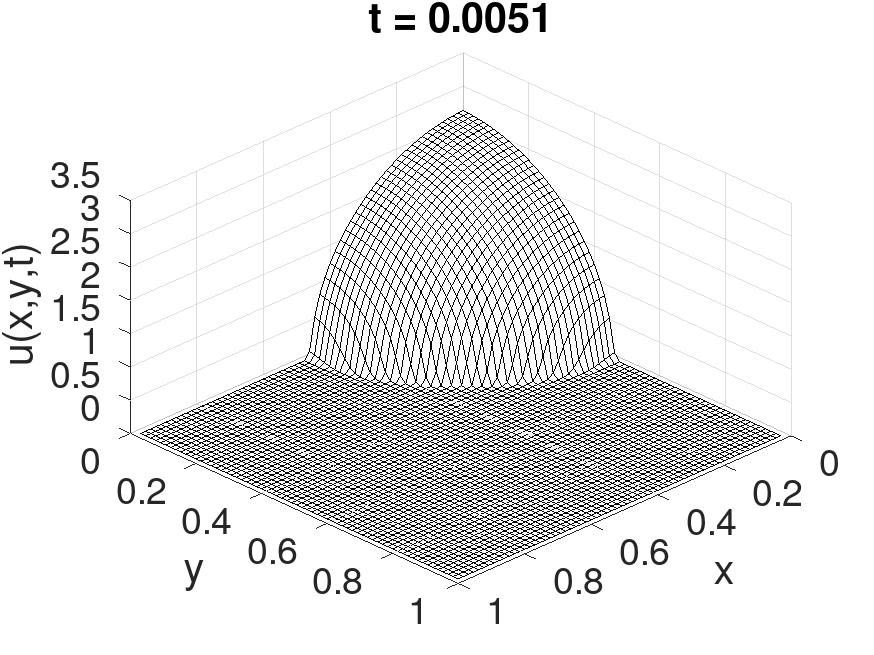}
\caption{Solution of the nonlinear heat equation at time
$t=0.0025$ (left plot) and $t=0.0051$ (right plot)}
\label{f:heat_sol}
\end{figure}

%
% ~/wrk/doc/matfun/matfun_split/resid_rst/nonlin/t5_heat2D$ ls -lath test_*.m
% -rw-r--r-- 1 botchev botchev 2.3K Mar 18 16:14 test_BE_adapt.m
% -rw------- 1 botchev botchev 2.2K Mar 18 11:55 test_LIM_adapt.m
% -rw-r--r-- 1 botchev botchev 2.7K Mar 18 11:39 test_exp_adapt_rst1.m
%
We now consider another test problem, an initial-boundary value problem
for nonlinear heat conduction equation
\begin{equation}
\label{heat}
\frac{\partial}{\partial t}u(x,y,t) = \nabla\cdot (k(u)\nabla u(x,y,t)),
\quad k(u) = u^2,
\end{equation}
where initial and boundary conditions are taken according to the exact solution
$u_{\mathrm{ex}}(x, y, t)$, chosen in this test as
\begin{equation*}
\label{uex} 
u_{\mathrm{ex}}(x, y, t) = \dfrac{1}{\sqrt[3]{t}}
\sqrt{\dfrac16 \max\left\{ 
0; \left(1.3 - \frac{x^2 +y^2}{\sqrt[3]{t}}\right) 
\right\} }.
\end{equation*}
The exact solution is shown, for two time moments, in Figure~\ref{f:heat_sol}.

The standard second order finite central difference approximation
on uniform Cartesian grid
reduces this initial-boundary value problem to an IVP
\begin{equation}
\label{heat_ivp}
y'(t) = -A(y(t))\, y(t) + g(t),\quad y(0) = y^0,
\end{equation}
where the sought after vector function $y(t)$ contains the grid values 
of the unknown
function~$u(x,y,t)$, $A(y)\in\Rr^{N\times N}$ is a symmetric positive 
semidefinite matrix approximating the terms $\nabla\cdot (k(u)\nabla u(x,y,t))$
in~\eqref{heat} and $g(t)$ is a known source vector function
incorporating the imposed boundary conditions.
This test problem is also considered in our earlier 
work~\cite{BotchevZhukov2023,BotchevZhukov2024}. 
     
To solve a problem class where problem~\eqref{heat} belongs to, 
in~\cite{BotchevZhukov2024} 
we proposed two stabilized explicit time integration schemes with
an adaptive time step size selection.  One of the schemes is the LIM
(local iteration modified) method, which is based on specially chosen
Chebyshev polynomials.  Assuming for simplicity of presentation that the matrix 
$A$ in~\eqref{heat_ivp} is ``frozen'' at some $\tilde{y}$, i.e.,
$A$ is set to $\tilde{A}=A(\tilde{y})$, and $g(t)\equiv 0$, 
consider the implicit Euler scheme for IVP~\eqref{heat_ivp},
\begin{equation}
\label{BE_lin}
 y^{n+1} = (I+\tau \tilde{A})^{-1} y^n , 
\quad n=0,1,2,\dots,
\end{equation}
where $\tau>0$ is the time step size.  In the LIM scheme, the
action of the inverse matrix $(I+\tau\tilde{A})^{-1}$ in~\eqref{BE_lin}
is replaced by $2p-1$ explicit Chebyshev iterative steps.
This means that the LIM solution
at the next time step~$n+1$ can formally be written as
\begin{equation*}
\label{LIMstep}  
y^{n+1} = \left( I-F_{p}^{2} (\tilde{A})\right)\, 
\left(I+\tau\tilde{A}\right)^{-1} y^n ,
\end{equation*}
where $F_p$ is a scaled Chebyshev polynomial
$$
F_p (\lambda) = \frac{1}{H_p (-1/\tau )} H_p (\lambda ), 
\quad  
H_{p} (\lambda )=T_p \left(z_1 - (z_{1} +1)\frac{\lambda}{\lambmax}\right),
$$
normalized by requirements $F_{p} (-1/\tau )=1$ and $F_{p} (0)=0$.
Here $\lambmax$ is an estimate of the largest eigenvalue of $\tilde{A}$.
The polynomial $F_p$ is constructed to be the least deviating 
polynomial for $\lambda\in[\lambda_0 \, ; \, \lambmax ]$, where
$$
\lambda_0 = \lambmax  \frac{ z -1 }{z +1} \in [-1/\tau ,0],
\quad
z =\cos \frac{\pi}{2p}.
$$
It can also be shown that 
$|F_p (\lambda )| \leqs 1$  for $\lambda\in [0 \, ; \, \lambmax ]$
and the LIM scheme is stable, i.e.,
$$
\left\|
\left( I-F_{p}^{2} (\tilde{A})\right)\, 
\left(I+\tau\tilde{A}\right)^{-1}\right\| \leqs 1,
$$      
provided that $p$ is taken sufficiently large, i.e.,
\begin{equation}
\label{p}
p = \left\lceil\dfrac{\pi/4}{\pi/2-\arctan\sqrt{\tau\lambmax}}\right\rceil.
\end{equation}
Here $\lceil x \rceil$ denotes the ceiling function (which is 
the smallest integer not less than~$x$).
Per time step the LIM scheme requires $2p-1$ matrix-vector products
with $\tilde{A}$ and it can be checked
that $1/(\pi/2-\arctan\sqrt{\tau\lambmax}) \leqs \sqrt{\tau\lambmax + 1}$.
Hence, the number of matrix-vector products in the LIM scheme,
determined by~\eqref{p}, is roughly proportional to $\sqrt{\tau\lambmax}$.
More details on the LIM scheme can be found in~\cite{Zhukov2011}.

To treat the nonlinearity in~\eqref{heat_ivp} with the LIM scheme, 
in~\cite{BotchevZhukov2024}, we applied a simple iterative scheme 
adopted from~\cite[App.~1, Ch.~2.11]{TikhonovSamarskii_UMF}.  
At each iteration $m=0,1,2,\dots$, where $y^{(0)}:=y^n$ and 
we hope to have $y^{(m)}\rightarrow y^{n+1}$, we solve
\begin{equation}
\label{BEiter}
y^{(m+1)} = ( I+\tau A(y^{(m)}) )^{-1}\left[ y^n + \tau g^{n+1} \right].
\end{equation}
The action of the inverse $( I+\tau A(y^{(m)}) )^{-1}$ here is replaced
by the LIM scheme step described above.

A time adaptive variant of this nonlinear LIM scheme, employing 
a predictor-corrector error estimation procedure, was tested and shown 
to work successfully in~\cite{BotchevZhukov2024}. 
We now examine how switching in this adaptive time integration
scheme to an implicit Euler step will influence its performance.
More specifically, to switch to implicit Euler, 
we just implement the nonlinear iteration~\eqref{BEiter} 
as it reads, i.e., we solve the linear system with the matrix
$I+\tau A(y^{(m)})$ by a direct sparse solver. 
For reference purposes, we also include in the comparisons 
the adaptive exponential Euler (EE) scheme proposed 
in~\cite{BotchevZhukov2024}.
All the test runs are carried out with the accuracy tolerance~$10^{-2}$,
for time interval $t\in[0.0001,0.0051]$.
The results are presented in Table~\ref{t:heat}.
The reported error values are relative grid $L_2$ error norms 
${\|y^n - y(T)\|}/{\|y(T)\|}$, where $y(T)$ is the exact solution
vector at the final time~$T=0.0051$ and $y^n$ is the numerical solution
at the last time step.

\begin{table}
\caption{Nonlinear heat equation problem.
Number of time steps and relative error norms
for the adaptive local iteration modified (LIM), 
implicit Euler and exponential Euler (EE) time integration schemes.
The shown stiffness values are estimates 
$\max_n\tau_n\rho(A(y^n))$, with $\rho(A)$ denoting the spectral radius
of $A$ and maximum taken among all time steps~$n$ for the implicit Euler scheme.}
\label{t:heat}
\begin{center}
\addtolength{\tabcolsep}{2ex}
\centerline{\begin{tabular}{cccc}
\hline\hline
grid size, & LIM scheme  & implicit Euler & EE scheme \\[-1ex]
stiffness $\displaystyle\max_n \tau_n\rho(A(y^n))$
& \# time steps, \texttt{error} 
& \# time steps, \texttt{error} 
& \# time steps, \texttt{error} 
\\\hline
$128\times 128$, 120 &
  671, {\tt6.88e-03} &   618, {\tt6.81e-03} &  370, {\tt7.03e-03} 
\\
$256\times 256$, 490 &
 1064, {\tt2.87e-03} &  1023, {\tt3.42e-03} &  663, {\tt4.75e-03}
\\
$512\times 512$, 1950 &
 1738, {\tt1.54e-03} &  1706, {\tt2.31e-03} & 1158, {\tt4.60e-03}
\\
$1024\times 1024$, 4650 &
 2881, {\tt9.86e-04} &  2856, {\tt1.67e-03} & 1985, {\tt5.06e-03}  
\\\hline
\end{tabular}}
\end{center}
\end{table}

\begin{figure}
\begin{center}
\includegraphics[width=0.49\linewidth]{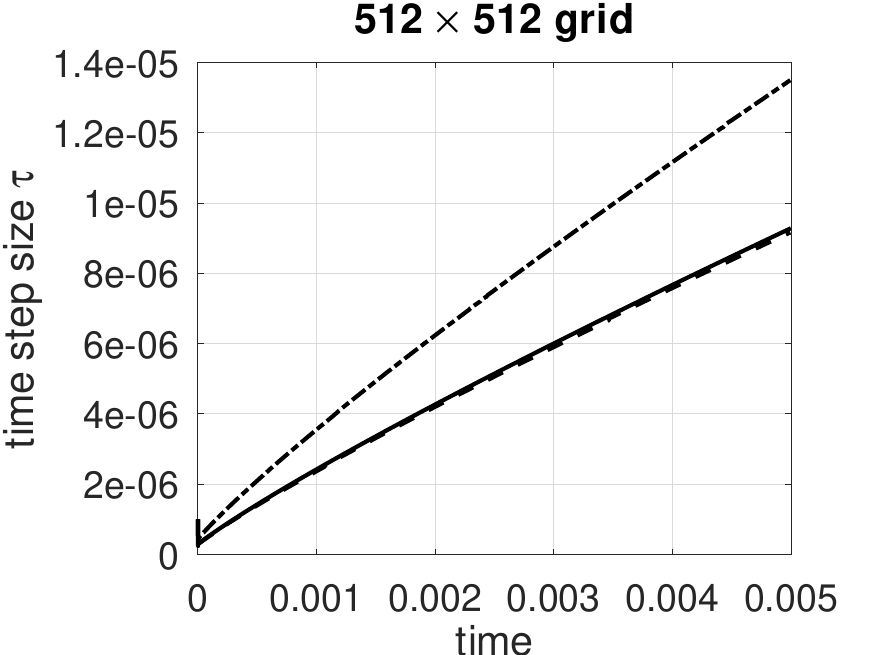}
\end{center}
\caption{Time step evolution in the implicit Euler scheme (solid line),
the LIM scheme (dashed line) and the EE scheme (dash-dotted line)
for the $512\times 512$ grid. The first two lines are hardly
distinguishable.}
\label{f:dt_heat}
\end{figure}

As we see in Table~\ref{t:heat}, 
for this test problem, employing the implicit Euler scheme leads to only
a slight 
%% %Correction to be hopefully made in the proofs: !!!
decrease in the number of %increase of the 
time steps needed to arrive at the final time~$T$.
Apparently, in this case the time step size is primarily restricted
by accuracy rather than by stability requirements.
Numerical tests with other tolerance values reveal that the adaptive
implicit Euler scheme is also slightly more robust than the adaptive LIM
scheme in terms of convergence of nonlinear iteration~\eqref{BEiter}.  
Time step size evolution for the three tested scheme is shown
in Figure~\ref{f:dt_heat}.  As we see, there is hardly any
difference in time step sizes histories of the LIM and implicit Euler
schemes.

%\textcolor{blue}{
Last but not least, we should emphasize that the results 
obtained in this section
may depend on the specific nonlinearity type of~$k(u)$ in~\eqref{heat} 
as well as on the way the nonlinearity is treated (i.e., in this case Richardson
iterations~\eqref{BEiter}).
%}

\section{Conclusions}
\label{s:concl}
The following conclusions can be drawn from the results presented
in this paper.
\begin{enumerate}
\item 
There is a clear gap in performance between the adaptive implicit trapezoidal
rule scheme and the explicit stabilized scheme based on diagonally
preconditioned BiCGstab iterations.  The latter scheme needs hundreds to
thousands matrix-vector multiplications per time step to have the same time 
step size as in the implicit scheme.  However, this performance gap turns out to
be much smaller than can be expected: 
%\textcolor{blue}{
for
%} 
the stiffness estimate $\tau\rho(M^{-1}A)$ 
%\textcolor{blue}{
reaching
%} 
the value $\approx 10^8$, around~$1700$  
matrix-vector multiplications are needed in the explicit stabilized scheme.

\item
There is, thus, a clear need for bridging the indicated gap,
i.e., for developing time integration schemes which are cheaper per time step
than 
%\textcolor{blue}{
regular
%} 
implicit schemes, yet sufficiently stable and parallelizable.

\item
The BiCGstab iterations employed in the stabilized explicit scheme
appear to be an order of magnitude more efficient than Chebyshev
polynomial iterations, which are rather often used in
stabilized explicit time stepping~\cite{RKC2004,Zhukov2011,Lebedev98}.  
This is apparently caused by the ability of Krylov subspace methods 
to adapt to the initial vector and to the discrete structure of the 
spectrum of the matrices involved.

\item
For the nonlinear heat equation problem, switching from the stabilized 
explicit local iteration modified (LIM) scheme to the implicit Euler
scheme appears not to 
%\textcolor{blue}{
yield
%} 
a significant increase in the time step size.  Apparently, in this problem 
class the time step size is restricted primarily by the accuracy rather than 
by stability requirements.

\item
%\textcolor{blue}{
The last conclusion should be accepted as holding 
for the specific nonlinearity type in~\eqref{heat} 
and the way the nonlinearity is treated, see~\eqref{BEiter}.
%}

\end{enumerate}

\section*{Acknowledgments and statements}
%\textcolor{blue}{
The author thanks the anonymous reviewer for useful
suggestions which helped to improve the paper and
Dmitry Zheltkov for fruitful discussions.

The author declares that no neural network techniques have been used 
for producing results reported in this work and writing the paper. 
%}

\section*{Funding}
This work is supported by ongoing research funding within the State assignment 
framework of the Keldysh Institute of Applied Mathematics of Russian Academy of 
Sciences.

\section*{Conflict of interest}
The author of this work declares that he has no conflicts of
interest.